\documentclass[12pt,english]{article}
\usepackage[cp1251]{inputenc}
\usepackage{amssymb,amsmath,babel,amsthm,graphicx,cite}
\usepackage{color}
\usepackage{epstopdf}

\usepackage{tikz}
\usetikzlibrary{graphs}

\newcounter{parag}

\newtheorem*{theorem*}{Theorem}

\theoremstyle{definition}

\newcounter{cl}

\begin{document}
	
	\begin{center}
		\large \textbf{On knot groups acting on trees}
		\normalsize \vspace{5mm}
		
		\textsc{F.~A.~Dudkin$^1$, A.~S.~Mamontov$^1$}
		
        \footnotetext[1]{The work was supported by the program of fundamental scientific researches of the SB RAS No I.1.1., project No 0314-2016-0001}
        
		\end{center}
\vspace{5mm}

{\bf Abstract:} A finitely generated group $G$ acting on a tree with infinite cyclic edge and vertex stabilizers is called a generalized Baumslag--Solitar group  {\it (GBS group)}.
We prove that a 1-knot group $G$ is GBS group iff $G$ is a torus-knot group and describe all n-knot GBS groups for $n\geqslant 3$.


\section{Introduction}
	
An {\it $n$-knot group } is the fundamental group $\pi_1 (S^{n+2}-K^n,a)$ for an $n$-knot $K^n$ in $n+2$-dimenstional sphere $S^{n+2}$. Starting from a knot it is possible to construct Wirtinger
presentation for its group with relations of the form $x_i^w=x^j$, where $x_i, x_j$ are letters, and $w$ is a word \cite{Kaw}. Now, let $G$ be a group represented by a set of generators and relations. When it is a knot group? Treating this just as a question of how to transform relations to a desired form is unfruitful. So to attempt the question we involve some known properties of knot groups
and restrict the class of groups.

A finitely generated group $G$ acting on a tree with infinite cyclic edge and vertex stabilizers is called a generalized Baumslag--Solitar group  {\it (GBS group)} \cite{forester2002}. By the Bass-Serre theorem, $G$ is representable as $\pi_1(\mathbb{B})$, the fundamental group of a graph of groups $\mathbb{B}$ (see \cite{Serr}) with infinite cyclic edge and vertex groups.

GBS groups are important examples of JSJ decompositions. JSJ decompositions appeared first in 3-dimensional topology with the theory of the characteristic
submanifold by Jaco-Shalen and Johannson. These topological ideas were carried over to group theory by Kropholler for some Poincar\'{e} duality groups of dimension at least 3, and by Sela for torsion-free hyperbolic groups. In this group-theoretical context, one has a finitely generated group G and a class
of subgroups $\mathcal{A}$ (such as cyclic groups, abelian groups, etc.), and one tries to understand splittings (i.e. graph of groups decompositions) of $G$ over groups in $\mathcal{A}$ (see \cite{GuLe}).

Given a $GBS$ group $G$, we can present the corresponding graph of groups by a labeled graph $\mathbb{A}=(A,\lambda)$, where $A$ is a finite connected graph (with endpoint functions $\partial_0,\partial_1\colon E(A)\to V(A)$) and $\lambda\colon E(A)\to \mathbb{Z}\setminus\{0\}$ labels the edges of $A$. The label $\lambda_e$ of an edge $e$ with the source vertex $v$ defines an embedding $\alpha_e\colon e\to v^{\lambda_e}$ of the cyclic edge group $\langle e \rangle$ into the cyclic vertex group $\langle v \rangle$ (for more details see \cite{ClayForII})

{\it The fundamental group} $\pi_1(\mathbb{A})$ of a {\it labeled graph} $\mathbb{A}=(A, \lambda)$ is given by generators and defining
relations. Denote by $\overline{A}$ the graph obtained from $A$ by identifying $e$ with $\overline{e}$. A maximal subtree $T$ of the graph $\overline{A}$ defines the following presentation of the group $\pi_1(\mathbb{A})$

$$\left\langle
\begin{array}{lcl}
g_v, v\in V(\overline{A}), &\|&  g_{\partial_0(e)}^{\lambda(e)}=g_{\partial_1(e)}^{\lambda(\overline{e})}, e\in E(T),\\
t_e, e\in E(\overline{A})\setminus E(T) &\|&t_e^{-1}g_{\partial_0(e)}^{\lambda(e)}t_e=g_{\partial_1(e)}^{\lambda(\overline{e})}, e\in E(\overline{A})\setminus E(T)
\end{array}
\right\rangle$$
Generatos of first (second) type are called vertex (edge) elements. For different maximal subtrees, corresponding presentations define isomorphic groups.

It is sometimes useful to regard a GBS-group as a group obtained as follows: start with the group $\mathbb{Z}$,
perform consecutive amalgamated products in accordance with the labels on the maximal subtree; finally,
apply several times the construction of the HNN-extension (the number of times is equal to the number
of the edges outside the maximal tree). In this approach, the standard theory of amalgamated products and HNN-extensions is applicable to the full extent. In particular, GBS-groups admit a normal form of an element and have no torsion.

If two labeled graphs $\mathbb{A}$ and $\mathbb{B}$ define isomorphic $GBS$ groups $\pi_1(\mathbb{A})\cong\pi_1(\mathbb{B})$ and $\pi_1(\mathbb{A})$ is not isomorphic to $\mathbb{Z}, \mathbb{Z}^2$ or Klein bottle group then there exists a finite sequence of {\it expansion} and {\it collapse} (see fig.1) moves connecting $\mathbb{A}$ and $\mathbb{B}$ \cite{forester2002}. A labeled graph is called {\it reduced} if it admits no collapse move (equivalently, the labeled graph contains no edges with distinct endpoints and labels $\pm 1$).

\begin{figure}[h]
\centering{\includegraphics[width=80mm]{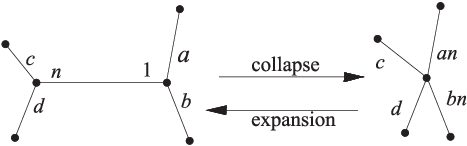}}
\caption{Expansion and collapse moves}
\end{figure}

An element $g$ from a GBS group $G$ is called elliptic if $g$ is conjugated with $a^k$, for some vertex generator $a$. The set consisting of all non-trivial elliptic elements is stable under conjugation, its elements have infinite order, and any two such elements are commensurable. These properties yield a homomorphism $\Delta$ from $G$ to the multiplicative group of non-zero rationals $\mathbb{Q}^\ast$, defined as follows.

Given $g\in G$, choose any non-trivial elliptic element $a$. There is a relation $g\cdot a^p\cdot g=a^q$, with $p, q$ non-zero, and define $\Delta(g)=\frac{p}{q}$. As pointed out in \cite{KropII} this definition is independent of the choices made ($a$ and the relation), and defines a {\it modular homomorphism}.

For different primes $p$ and $q$ let $T(p,q) = \langle x,y | x^p=y^q \rangle$ be a torus-knot group. It is easy to see that $T(p,q)$ is a $GBS$ group for any non-zero integers $p$ and $q$ (corresponding labeled graph has one edge with two different endpoints and labels $p$ and $q$).

A group is said to be {\it Hopfian} if any homomorphism of the group onto itself has trivial kernel, i.e.
is an automorphism. Baumslag and Solitar \cite{BaSo} came up with a series of examples of two-generator
one-relator non-Hopfian groups. In particular, such are the
Baumslag--Solitar groups
$$BS(p,q)=\langle x,y |xy^{p}x^{-1}=y^{q} \rangle$$
where $p$ and $q$ are coprime integers, $p,q\neq 1$.

If a labeled graph $\mathbb{B}$ consists of one vertex and two inverse loops with labels $p$ and $q$, then
$\pi_1(\mathbb{B})\cong BS(p,q)$. Therefore, every Baumslag--Solitar groups is a generalized Baumslag--Solitar
group.

Main results of our work are listed below.
	
{\bf Theorem 1.} {\it Let $G$ be a GBS-group. Then $G$ is 1-knot group if and only if $G \simeq T(p,q)$.}

{\bf Theorem 2.} {\it Let $G$ be a GBS-group, $G\not\cong\mathbb{Z}$. Then $G$ is $n$-knot group for $n \geq 3$ if and only if $G$ is a homomorphic image of either $BS(m,m+1)$, where $m\geq 1$, or $T(p,q)$.}

All homomorphic images of $BS(m,m+1)$, where $m\geq 1$, and $T(p,q)$ described in terms of labeled graphs in lemmas 4 and 5 (see the proof of theorem 2).

M.Kervaire obtained the following necessary conditions for a group to be a knot-group.

{\bf Statement 1 \cite[14.1.1]{Kaw}}. {\it Let $\pi$ be $n$-knot group for $n \geq 1$. Then

1. $\pi$ is finitely generated

2. $\pi / \pi ' \simeq \mathbb{Z}$

3. $H_2(\pi)=0$

4. $\pi$ is a normal closure of a single element }

For $n \geq 3$ these conditions are also sufficient.

{\bf Statement 2 \cite[14.1.2]{Kaw}}. {\it If a group $G$ satisfies conditions 1-4 of lemma 1, then $G$ is $n$-knot group for $n\geq 3$.}

So our goal is to determine when these conditions are fulfilled in a GBS group.

\section{Preliminary results}

Given graph $A$, denote the number of edges out of maximal subtree by $\beta_1(A)$. Normal closure of an element $a$ in a group $G$ is denoted as  $\langle \langle a \rangle \rangle _G$.

{\bf Lemma 1.} {\it Let $G$ be a GBS group such that $G/G' \simeq \mathbb{Z}$, then $\beta_1(A) \leqslant 1$.}

Proof.  From a representation it is clear that $G/G' \simeq \mathbb{Z}^{\beta_1(A)}\times H$, where $H$ is a subgroup of $G/G'$ generated by images of vertex elements. Hence, $\beta_1(A) \leqslant 1$.
Lemma is proved.

{\bf Lemma 2.} {\it Let $G=\pi_1(\mathbb{A})$ be a $GBS$-group such that $G= \langle \langle a \rangle \rangle _G$, and $\mathbb{A}$ be a labeled tree. Then $\mathbb{A}$ is a labeled segment.}

\begin{figure}[h]
\centering{\includegraphics{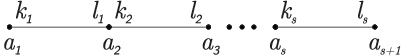}}
\caption{Labeled segment.}
\end{figure}

Proof. Assume the contrary. Then  $\mathbb{A}$ has a trident subgraph.

\begin{figure}[h]
\centering{\includegraphics{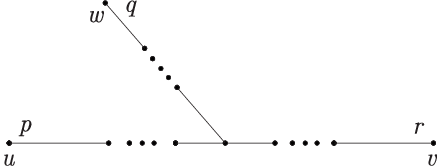}}
\caption{Labeled trident.}
\end{figure}

Let $N$ be a normal closure of all vertex elements except end points $u,v,w$, then $G/N \simeq \mathbb{Z}_p * \mathbb{Z}_q * \mathbb{Z}_r$.

If $p,q,r$ are pairwise comprime, then for some normal
subgroup $N_1$ containing $N$ and certain powers of $u,v,w$ we have $G/N_1 \simeq \mathbb{Z}_{p_1} * \mathbb{Z}_{q_1} * \mathbb{Z}_{r_1}$, where $p_1,q_1,r_1$ are different primes. By \cite{How} $G/N_1$ is not a normal closure of a single element, so neither is $G$, a contradiction.

Now let $(p,q)=d \not = 1$.
GBS group $G= \langle \langle a \rangle \rangle _G$ has no torsion, therefore $G/G' \simeq \mathbb{Z}$. Thus we have $(G/G')/N \simeq Z_p \times Z_q \times Z_r$. With $d \not =1$ dividing $p$ and $q$ the later group cannot be a homomorphic image of $\mathbb{Z}$.

Using ideas of plateau from \cite{Levitt2015}, we prove.

{\bf Lemma 3.} {\it Keep notations of lemma 2. Then $l_i \bot k_j$ for all $i,j$.}

Proof. Suppose first that $j \leqslant i$.

Choose a pair $k_j, l_i$ so that $j \leqslant i$, $(k_j,l_i)=d>1$, $j$ is minimal and $i$ is maximal with these conditions. Let $p$ be a prime divisor of $d$. Then $k_x \perp p$ for $x < j$ and $l_y \perp p$ for $y > i$. Let $r$ be a maximal index such that $l_{r-1} \not \perp p$, $2 \leq r \leq j$. If there is no such index, then let $r=1$. In a similar way, choose minimal $f$, $f \geq i+1$ such that $p | k_f$.

\begin{figure}[h]
\centering{\includegraphics{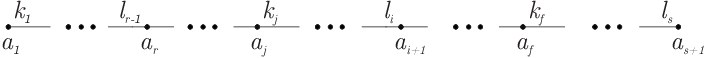}}
\caption{Labeled segment, $j \leqslant i$.}
\end{figure}

Let $N$ be a normal closure of the following elements

$$N = \langle \langle a_1, \ldots a_{r-1}, a_{j+1}, \ldots a_i, a_{f+1} \ldots a_{s+1},a_r^p,\ldots a_j^p,a_{i+1}^p,\ldots , a_f^p \rangle \rangle$$

Then corresponding factor group $G_1=G / N$ is isomorphic to $\mathbb{Z}_p * \mathbb{Z}_p$, since relations

$$a_m^{k_{m+1}} = a_{m+1}^{l_{m+1}},$$

$$a_m^p = a_{m+1} ^p=1$$

and $p \perp k_{m+1}, p \perp l_{m+1}$ imply $a_m=a_{m+1}$, $m=r,\ldots j,i+1,\ldots , f$.

Hence $G_1 / G_1' \simeq \mathbb{Z}_p \times \mathbb{Z}_p$, which condtradicts $G/G' \simeq \mathbb{Z}$.

2. Now assume $i<j$ and $(k_j,l_i)=d>1$.

\begin{figure}[h]
	\centering{\includegraphics{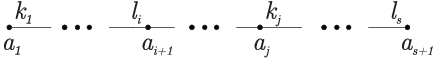}}
	\caption{Labeled segment, $j>i$.}
\end{figure}

We may also assume that $j-i$ is minimal with this property. Then using 1 we obtain $k_r \perp l_f$, $r=i+1,\ldots j-1, f=i,\ldots , j-1$ and $k_j \perp l_{i+1}, \ldots l_{j-1}$. In particular, for a prime divisor $p$ of $d$ we have $p \perp k_{i+1}, \ldots k_{j-1}$ and $p \perp l_{i+1}, \ldots , l_{j-1}$. Let

$$N= \langle \langle a_2,\ldots , a_i, a_{j+1}, \ldots , a_s, a_{i+1}^p,\ldots , a_j^p \rangle \rangle.$$

Then as above we obtain

$$G/N \simeq \mathbb{Z}_{k_1} * \mathbb{Z}_p * \mathbb{Z}_{l_{s+1}}$$

Having $k_1 \perp l_i$ and $k_j \perp l_{s+1}$ we may assume that $k_1,p,l_{s+1}$ are different primes. By our assumption $G$, and hence $G/N$, is a normal closure of a single element, while
$ \mathbb{Z}_{k_1} * \mathbb{Z}_p * \mathbb{Z}_{l_{s+1}}$ for different primes  $k_1,p,l_{s+1}$ is not \cite{How}. This contradiction proves the lemma.

{\bf Lemma 4.} {\it Assume $G = \pi _1 (\mathbb{A})$,  $\mathbb{A}$ is a labeled segment, and $l_i \perp k_j$ for all $i,j$. Then $G$ is a homomorphic image of some torus-knot group. }

Proof.

We prove by induction that $a_2 , \ldots a_{s-1}$ may be excluded from generators.

To check basis consider a group with generators $\{ a_1,a_2,a_3 \}$ and relations $a_1^{k_1}=a_2^{l_1}$, $a_2^{k_2}=a_3^{l_2}$. Choose $\alpha$ and $\beta$ such that $\alpha l_1+\beta k_2 =1$. Then $a_2=a_1^{\alpha k_1} a_3^{\beta l_2 }$.

For the induction step consider a group with generators  $\{ a_1,a_2, \ldots, a_{s+1} \}$. By induction, $a_i = w_i (a_1,a_s)$, where $w_i$ is some word on letters $a_1,a_s$, and $i=2,\ldots s-1$.  Choose $\alpha$ and $\beta$ so that $\alpha l_1\ldots l_{s-1}+\beta k_s =1$. Then $a_s=a_s^{\alpha l_1 \ldots l_{s-1}} a_s^{\beta k_s} = a_1^{\alpha k_1 \ldots k_{s-1}} a_{s+1}^{\beta l_s}$. So a group is generated by $a_1$ and $a_{s+1}$, which satisfy a relation $a_1^k=a_{s+1}^l$, with $k=k_1\ldots l_s \perp l=k_1\ldots l_s$. Therefore $G$ is a homomorphic image of some $\pi_1 (T_{k,l})$. The
lemma is proved.

{\bf Lemma 5.} {\it  If $G=\pi_1(\mathbb{A})$, $G=\langle \langle a \rangle \rangle_G$, and $\beta_1(A)=1$ then $\mathbb{A}$ is a labeled cycle with coprime labels $k_i \perp l_j$ for all $i,j$, $|{\displaystyle \prod_{i=1}^{s} k_i} - {\displaystyle \prod_{i=1}^{s} l_i}|=1$ (see fig.6) and $G= \langle \langle t \rangle \rangle_G$.}

\begin{figure}[h]
	\centering{\includegraphics{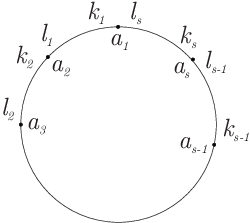}}
	\caption{Labeled cycle}
\end{figure}

Proof.
If $\mathbb{A}$ is no a cycle, then there is a pendant vertex $u$. A graph  $\mathbb{A}$ is reduced, and so the label $\lambda$ near $u$ is not $\pm 1$. Let $N$ be a normal closure of all vertices, except $u$.
Then $G / G'N \simeq \langle t \rangle \times \langle u \rangle  \simeq \mathbb{Z} \times \mathbb{Z}_{\lambda}$, a contradiction.

Let $d=(k_i,l_j)$ and $p$ be a prime divisor of $d$. In case $j <i$ let $N$ be a normal closure of $a_1,\ldots a_j,a_{j+1}^p,\ldots ,a_i^p,a_{i+1},\ldots , a_s$. Then
$G / G'N \simeq \langle t \rangle \times \langle a_{j+1} \rangle  \simeq \mathbb{Z} \times \mathbb{Z}_p$, as in lemma 3, a contradiction. The case $j \geq i$ is considered in a similar way.

Let $k= {\displaystyle \prod_{i=1}^{s} k_i}$, $l= {\displaystyle \prod_{i=1}^{s} l_i}$.
If $k=1$, then $s=1$, because  $\mathbb{A}$ is reduced. So $G \simeq \pi_1(\mathbb{A}) \simeq BS(1,l)$, $G /G' \simeq \mathbb{Z} \times \mathbb{Z}_{l-1}$. It follows that $l=2$, and $BS(1,2) = \langle \langle t \rangle \rangle$ satisfies the conclusion of the lemma.

If $k\not =1$ and $l\not =1$. Note that in $G$ we have an equality $t^{-1}a_i^kt=a_i^l$ for  $1 \leq i \leq s$, and $l \perp k$.

A map $\phi _{[a_1,t]}: H=BS(k,l)=\langle a,r| r^{-1}a^kr=a^l \rangle \rightarrow G$

$$ \phi _{[a_1,t]}: \begin{cases}
a \rightarrow a_1\\
r \rightarrow t[a_1,t] \end{cases}$$

is an embedding $BS(k,l) \hookrightarrow G$ \cite{DBS}, hence $|a_1|_{G/G'}=|a|_{H/H'}=|k-l|$. Therefore $|k-l|=1$. Moreover there is a subgroup $\langle t \rangle \times \langle a_1 \rangle$ in $G/G'$ isomorphic to $\mathbb{Z} \times \mathbb{Z}_{|k-l|}$. Note that $t^{-1}a_i^kta_i^{-k} \in \langle \langle t \rangle \rangle$, hence $a_i^{l-k} = a_i^{\pm 1} \in \langle \langle t \rangle \rangle$ for all $i$, and so $G = \langle \langle t \rangle \rangle$. The lemma is proved.

\section{Proof of theorems}

{\bf Proof of theorem 1}

Assume that $G$ is both a GBS-group and a 1-knot group.

1. By statement 1 $G/G' \simeq \mathbb{Z}$.

2. By lemma 1 $G\simeq \pi_1(\mathbb{A})$ and either $A$ is a tree or $\beta_1(A)=1$.

3. If $A$ is a tree, then $\Delta(G)=\{1\}$ and by \cite[Statement 2.5]{Levitt2007} $Z(G)=Z$.

4. By \cite[Corollary 6.3.6]{Kaw} $G\simeq T(p,q)$.

5. If $\beta_1(A)=1$, then by \cite[Theorem 6.3.9]{Kaw} $G=\pi_1 (\mathbb{A})$ is residually finite. By \cite[Corollary 7.7]{Levitt2015q} a GBS-group is residually finite iff either $G=BS(1,n)$
or $\Delta(G)=\{\pm 1\}$.

6. If $G=BS(1,n)$ then by Lemma 5 $n=2$ and $G=BS(1,2)$, and by \cite{Shalen} $|m|=|n|$, a contradiction.

7. If $\Delta(G)=\{\pm 1\}$, then there exist $a,g \in G$ such that $gag=a^{-1}$, and hence $G/G'$ has a torsion, a contradiction with $G/G' \simeq \mathbb{Z}$.

8. If $\Delta(G)=\{1\}$, then $Z(G)=\mathbb{Z}$ \cite[Statement 2.5]{Levitt2007}, and $G$ is isomorphic to $T(p,q)$ by \cite[Corollary 6.3.6]{Kaw}.

The Theorem 1 is proved.

{\bf Proof of theorem 2}

I. Assume that $G$ is a GBS-group, $G \not \simeq \mathbb{Z}$ and $G$ is $n$-knot group, $n\geq 3$.

1. By Statement 1 $G/G' \simeq \mathbb{Z}$, $G=\langle \langle a \rangle \rangle$.

2. By Lemmas 1-5 $G$ is a fundamental group of either some labeled segment (see fig.2) with $k_i \perp l_i$, or some labeled cycle (see fig.6)  $k_i \perp l_i$.

3. If $G$ is a fundamental group of a segment on fig.2, then by Lemma 4 $G$ is a homomorphic image of $T(k,l)$.

4. If $G$ is a fundamental group of a cycle on fig.6 then by \cite[Theorem 1.1]{Levitt2015q} $G$ is a homomorphic image of a Baumslag Solitar group $BS(k,k+1)$.

II. Assume that $G$ is a homomorphic image of rank 2 of either $BS(m,m+1)$ or $T(p,q)$.

5.  If $G$ is a homomorphic image of rank 2 of $BS(m,m+1)$, then by \cite[Theorem 1.1]{Levitt2015q}, using that $m \perp m+1$, we get $G=\pi_1(\mathbb{A})$, where $\mathbb{A}$ is a cycle, and
${\displaystyle \prod_{i=1}^{s} k_i} \perp {\displaystyle \prod_{i=1}^{s} l_i}$, i.e. $\mathbb{A}$ as in conclusion of Lemma 5. By \cite{Krop} we have $H_2(G)=0$. Then by Lemma 5 we obtain that properties 1-4 of Statement 1 hold, and hence by Lemma 2 $G$ is
a group of $n$-knot, $n\geq 3$.

6.  If $G$ is a homomorphic image of rank 2 of $T(p,q)$, then by \cite{Levitt2015} and \cite{Levitt2015q} $G=\pi_1(\mathbb{A})$, where $A$ is either a segment, or a circle or a lollipop. Moreover,
from \cite{Levitt2015} it follows that $A$ is a segment from theorem 1.1 \cite{Levitt2015}. Then $A$ satisfies the conclusion of Lemma 4. Hence $G$ satisfies the conditions 1-4 of Statement 1, and by Statement 2 $G$ is
$n$-knot group for $n\geq 3$.

The Theorem 2 is proved.

Acknowledgments. The authors are grateful to V.~G.~Bardakov and V.~A.~Churkin for valuable discussions.

\bibliographystyle{plain}

\end{document}